\documentclass[12pt]{amsart}
\usepackage{latexsym}
\usepackage[all]{xy}
\usepackage{amsfonts}
\usepackage{amsthm}
\usepackage{amsmath}
\usepackage{amssymb}
\usepackage{amscd}
 \numberwithin{equation}{section}

    \def\<{{\langle}}
    \def\>{{\rangle}}
    
    \def\eps{\varepsilon}

    \def\note#1{{}}

    \def\note#1{}
    \def\M{{\bf M}}

    \def\cC{{\mathfrak C}}
    
    \def\cD{{\mathfrak D}}
    
    \def\cG{{\mathcal G}}

    \def\roA{{\varrho^A}}

    \def\hom#1#2#3{{{\rm Hom}\sb{#1}(#2,#3)}}

    \def\Rend#1#2{{{\rm End}\sp{#1}(#2)}}
    
    \def\Rhom#1#2#3{{{\rm Hom}\sp{#1}(#2,#3)}}
     \def\Raut#1#2{{{\rm Aut}\sp{#1}(#2)}}
      \def\raut#1#2{{{\rm Aut}\sb{#1}(#2)}}
       \def\Twist#1#2{{{\rm Twist}(#1,#2)}}
       \def\Tors#1{{{\rm Tors}(#1)}}

    \def\beq{\begin{equation}}
    \def\eeq{\end{equation}}
    \def\DC{{\Delta_\cC}}
    \def \eC{{\eps_\cC}}

    \def\ot{{\otimes}}

    \def\roM{\varrho^{M}}

    \def\cten#1{\raise-.2cm\hbox{$\stackrel{\displaystyle\square}
{\scriptscriptstyle{#1}}$}}

    \def\tM{\widetilde{M}}

     \newcounter{zlist}

  \newcounter{blist}

\newtheorem{proposition}{Proposition}[section]
\newtheorem{lemma}[proposition]{Lemma}

\newtheorem{theorem}[proposition]{Theorem}
\theoremstyle{definition}
\newtheorem{definition}[proposition]{Definition}

\theoremstyle{remark}

\begin{document}
\title[{Descent cohomology and corings}]{Descent cohomology and corings}

   \author{Tomasz Brzezi\'nski}
   \address{ Department of Mathematics, University of Wales Swansea,
   Singleton Park, \newline\indent  Swansea SA2 8PP, U.K.}
   \email{ T.Brzezinski@swansea.ac.uk}

     \date{\today}
   \begin{abstract}
 A coring approach to non-Abelian descent cohomology of [P.\ Nuss \& M.\ Wambst, {\em Non-Abelian Hopf cohomology}, Preprint arXiv:math.KT/0511712, (2005)] is described and a definition of a Galois cohomology for partial group actions is proposed.
   \end{abstract}
   \maketitle

\section{Introduction}
\subsection{Motivation and aims}
In a recent paper \cite{NusWam:non} Nuss and Wambst have introduced a non-Abelian descent cohomology for Hopf modules and related it to classes of twisted forms of modules corresponding to a faithfully flat Hopf-Galois extension. Hopf modules can be understood as a special class of {\em entwined modules} and hence comodules of a coring \cite{Brz:str}. The aim of this note is to show how the descent cohomology 
introduced in \cite{NusWam:non} fits into recent developments in the descent theory for corings \cite{ElKGom:com}, \cite{Cae:Gal}, \cite{CaeDeG:com}. In particular we construct the zeroth and first descent cohomology sets for a coring with values in a comodule and relate it to isomorphism classes of module-twisted forms. We then use this general framework to propose a definition of a non-Abelian Galois cohomology for idempotent partial Galois actions on non-commutative rings introduced in \cite{CaeDeG:par}.
 
\subsection{Notation and conventions}\label{sec.not} We work over a commutative associative ring $k$ with unit. All algebras are associative, unital and over $k$. The identity map in a $k$-module $M$ is denoted by $M$. Given an algebra $A$ the coproduct in an $A$-coring $\cC$ is denoted by $\DC$ and the counit by $\eC$. 
A (fixed) coaction in a right $\cC$-comodule $M$ is denoted by $\roM$. $\Rhom\cC - -$, $\Rend\cC -$ and $\Raut \cC - $ denote the homomorphisms, endomorphisms and automorphisms of right $\cC$-comodules, respectively. $\Rend\cC M$ is a ring with the product given by the composition, $M$ is a left $\Rend\cC M$-module with the product given by  evaluation and $\Rhom\cC MN$ is a right $\Rend\cC M$-module with the product given by composition.

For an algebra $A$, $\cG(A)$ denotes the group of units in $A$, and for an $A$-coring $\cC$, $\cG(\cC)$ is the set of grouplike elements of $\cC$, i.e.\ elements $g\in G$ such that $\DC(g) = g\ot_A g$ and $\eC(g)=1$. $\cG(\cC)$ is a right $\cG(A)$-set with the action given by the conjugation $g\cdot u := u^{-1} g u$, for all $u\in \cG(A)$ and $g\in \cG(\cC)$.

More details on corings and comodules  can be found in  \cite{BrzWis:cor}.

\section{Descent and partial Galois cohomologies}
\subsection{Construction of descent cohomology sets}
Given an $A$-coring $\cC$ and a right $\cC$-comodule $M$ with fixed coaction $\roM: M\to M\ot_A\cC$, define the {\em zeroth descent cohomology group of $\cC$ with values in $M$} as the group of $\cC$-comodule automorphisms, i.e.
$$
D^0(\cC,M) := \Raut \cC M.
$$
Any isomorphism of right $\cC$-comodules, $f: \tM\to M$, induces an isomorphism of cohomology groups
$$
f^* : D^0(\cC,M)\to D^0(\cC,\tM), \qquad \alpha\mapsto f^{-1}\circ\alpha\circ f.
$$
The set $Z^1(\cC, M)$ of {\em descent 1-cocycles on $\cC$ with values in $M$} is defined as a set of all $\cC$-coactions $F:M\to M\ot_A\cC$. Since $M$ comes equipped with the right coaction $\roM$, $Z^1(\cC, M)$ is a pointed set with a distinguished point $\roM$.

Note that for the Sweedler coring $A\ot_BA$ associated to a ring extension $B\to A$, $Z^1(A\ot_B A, M)$ is the set of {\em non-commutative descent data} on $M$ (cf.\ \cite[Section~25.4]{BrzWis:cor}). This motivates the name {\em descent cocycles}.
\begin{lemma}\label{lemma.action}
Let $M,\tM$ be right $\cC$-comodules. Any right $A$-linear isomorphism $f:\tM\to M$ induces a bijection $f^*: Z^1(\cC,M)\to Z^1(\cC,\tM)$ defined by
$$
f^*(F) := (f^{-1}\ot_A\cC)\circ F\circ f.
$$
The operation $(-)^*$ maps the identity map  into the identity map and reverses the order of the composition, i.e., for all right $A$-module isomorphisms $f:\tM\to N$, $g:N\to M$, $(g\circ f)^* = f^*\circ g^*$. Furthermore, if $f: \tM\to M$ is an isomorphism of $\cC$-comodules, then $f^*$ is an isomorphism of pointed sets.
\end{lemma}
\begin{proof}
If $F: M\to M\ot_A\cC$ is a right coaction, and $f:\tM\to M$ is a right $A$-linear isomorphism, then
\begin{eqnarray*}
(f^*(F)\ot_A \cC)\circ f^*(F) &=& (f^{-1}\ot_A\cC\ot_A\cC)\circ (F\ot_A\cC)\circ (f\ot_A\cC)\circ(f^{-1}\ot_A\cC)\circ F\circ f\\
&=& (f^{-1}\ot_A\cC\ot_A\cC)\circ (M\ot_A\DC)\circ F\circ f\\
&=& (\tM\ot_A\DC)\circ (f^{-1}\ot_A\cC)\circ F\circ f = (\tM\ot_A\DC)\circ f^*(F),
\end{eqnarray*}
where the second equality follows by the coassociativity of $F$. This proves that $f^*(F)$ is a coassociative coaction. The counitality of $F$ immediately implies that also $f^*(F)$ is a counital map. Thus $f^*$ is a well-defined map, and it is obviously a bijection as stated. The proofs of remaining statements are straightforward.
\end{proof}

Lemma~\ref{lemma.action} immediately implies that 
there is a right action of the $A$-linear automorphism group $\raut A M$ on $Z^1(\cC,M)$ given by
$$
Z^1(\cC,M)\times \raut A M\to Z^1(\cC,M), \qquad (F,f)\mapsto f^*(F).
$$

\begin{definition}\label{def.desc}
The {\em first descent cohomology set of $\cC$ with values in the $\cC$-comodule $M$}  is
defined as the quotient of $Z^1(\cC,M)$ by the action of $\raut AM$
 and is denoted by $D^1(\cC,M)$.
\end{definition}

Since $Z^1(\cC,M)$ is a pointed set, so is $D^1(\cC,M)$  with the class of $\roM$ as a distinguished point.

\subsection{$M$-torsors}
Given a right $\cC$-comodule $M$, an {\em $M$-torsor} is a triple $(X, \varrho^X, \beta)$, where $X$ is a right $\cC$-comodule with the coaction $\varrho^X$ and $\beta:M\to X$ is an isomorphism of right $A$-modules. Two $M$-torsors $(X, \varrho^X,\beta)$ and $(Y, \varrho^Y,\gamma)$ are said to be {\em equivalent} if $(X,\varrho^X)$ and  $(Y,\varrho^Y)$ are isomorphic as comodules. 
Equivalence classes of $M$-torsors are denoted by $\Tors M$. $\Tors M$ is a pointed set with the class of the $M$-torsor $(M,\roM, M)$ as a distinguished point. Since $\Tors M$ is a set of isomorphism classes of comodules which are isomorphic to $M$ as modules, one obtains the following description of $\Tors M$.
\begin{proposition}\label{prop.tor}
For all $\cC$-comodules $M$, there is an isomorphism of pointed sets
$$
D^1(\cC,M)\simeq \Tors M.
$$
\end{proposition}
\begin{proof}
Explicitly, the isomorphism is constructed as follows. Given a coaction $F:M\to M\ot_A\cC$, consider an $M$-torsor $T(F):= (M,F,M)$. Conversely, given an $M$-torsor $(X, \varrho^X,\beta)$, define the coaction on $M$, 
$$D(X, \varrho^X,\beta) := \beta^*(\varrho^X) = (\beta^{-1}\ot_A\cC)\circ \varrho^X\circ\beta$$
(cf.\ Lemma~\ref{lemma.action}). 
\end{proof}

\subsection{Module-twisted forms}
For $k$-algebras $A$ and $B$, fix a $(B,A)$-bimodule $\Sigma$ and a right $B$-module $N$. A right $B$-module $P$ is called a {\em $\Sigma$-twisted form of $N$} in case there exists a right $A$-module isomorphism $\phi: P\ot_B\Sigma\to N\ot_B\Sigma$. $\Sigma$-twisted forms $(P,\phi)$, $(Q,\psi)$ of $N$ are said to be {\em equivalent} if  $P$ and $Q$ are isomorphic as right $B$-modules.
The equivalence classes of $\Sigma$-twisted forms of $N$ are denoted by $\Twist \Sigma N$. $\Twist \Sigma N$ is a pointed set with the class of $(N, N\ot_B \Sigma)$ as a distinguished point

\subsection{Descent cohomology and the Galois-comodule twisted forms.}
Recall that a right $\cC$-comodule $\Sigma$ is called a {\em Galois comodule} or $(\cC,\Sigma)$ is called a {\em Galois coring} in case $\Sigma$ is finitely generated and projective as an $A$-module and the evaluation map
$$
\Rhom \cC \Sigma \cC \ot_B \Sigma \to \cC, \qquad f\ot s\mapsto f(s),
$$
is an isomorphism of right $\cC$-comodules. Here and in the remainder of this subsection $B$ is the endomorphism ring $B:=\Rend \cC \Sigma$. 
\begin{theorem}\label{thm.desc.twist}
Let $\Sigma$ be a Galois right $\cC$-comodule that
is faithfully flat as a left $B$-module. Then, for all right $B$-modules $N$, there is an isomorphism of pointed sets
$$
D^1 (\cC, N\ot_B\Sigma)\simeq \Twist \Sigma N,
$$
where $N\ot_B\Sigma$ is a comodule with the induced coaction $N\ot_B\varrho^\Sigma$.
\end{theorem} 
\begin{proof}
First recall that by the Galois comodule structure theorem \cite[Theorem~3.2]{ElKGom:com} (cf.\ \cite[18.27]{BrzWis:cor}) the functors $\Rhom \cC \Sigma -$ and $-\ot_B\Sigma$ are inverse equivalences between the categories of right $\cC$-comodules and right $B$-modules. Take any right $B$-module $N$. Note that $\Tors {N\otimes_A\Sigma}$  are $\cC$-comodule isomorphism classes of comodules which are isomorphic to $N\ot_B\Sigma$ as right $A$-modules, while $\Twist \Sigma N$ are isomorphism classes of $B$-modules $P$ such that $P\ot_B\Sigma$ is isomorphic to $N\otimes_B\Sigma$ as a right $A$-module (and hence as a $\cC$-comodule, with induced coactions). Since $-\ot_B\Sigma$ is an equivalence, there is an isomorphism of pointed sets $\Tors {N\otimes_A\Sigma} \simeq \Twist \Sigma N$, and the assertion follows by Proposition~\ref{prop.tor}
\end{proof}
  
  \subsection{Beyond faithfully flat and finite Galois comodules} 
 Start with two algebras $A$ and $B$, an $A$-coring $\cC$ and  a $(B,A)$-bimodule $\Sigma$ with a left $B$-linear right $\cC$-coaction $\varrho^\Sigma: \Sigma\to \Sigma\ot_A\cC$. This defines the functor $-\otimes_B\Sigma: \M_B\to \M^\cC$, where $\M_B$ is the category of right $B$-modules and  $\M^\cC$ is the category of right $\cC$-comodules. If this functor is an equivalence, then the same reasoning as in the proof of Theorem~\ref{thm.desc.twist} yields
 an isomorphism of pointed sets
$$
D^1 (\cC, N\ot_B\Sigma)\simeq \Twist \Sigma N,
$$
where $N\ot_B\Sigma$ is a comodule with the induced coaction $N\ot_B\varrho^\Sigma$.

Let $B=\Rend \cC\Sigma$. For Galois comodules which are  finitely generated and projective as right $A$-modules but  {\em are not} faithfully flat as $B$-modules sufficient conditions for $-\ot_B\Sigma$ to be an equivalence are found in \cite[Theorem~4.6]{BohVer:Mor}.  In particular $
D^1 (\cC, N\ot_B\Sigma)\simeq \Twist \Sigma N,
$ for a {\em cleft bicomodule} $\Sigma$ for a right coring extension $\cD$ of $\cC$  provided $\cD$ has a grouplike element  (cf.\  \cite[Definition~5.1,~Corollary~5.5]{BohVer:Mor}). For Galois comodules in the sense of Wisbauer \cite{Wis:gal} (i.e.\ comodules $\Sigma$ such that the map of functors 
$\hom A \Sigma - \ot_{B}\Sigma\to -\ot_A\cC$, $\phi\ot s\mapsto (\phi\ot_A\cC)(\varrho^\Sigma(s))$ is a natural isomorphism, no finiteness assumption on $\Sigma_A$) the functor $-\ot_B\Sigma:\M_B\to \M^\cC$ is an equivalence if and only if the functor $-\ot_B\Sigma:\M_B\to \M_A$ is comonadic (cf.\  \cite[Section~1]{BohVer:Mor}). 

Even more generally,  $\Sigma$-twisted forms can be defined if the algebra $B$ is non-unital (e.g.\ $B$ can be an ideal in $\Rend \cC\Sigma$). If $B$ is {\em firm} in the sense that the product map $B\otimes_BB\to B$ is an isomorphism,  and $\Sigma$ is a firm $B$-module, then in the case of {\em infinite comatrix corings}  the  conditions for $-\ot_B\Sigma$ to be an equivalence are given in  \cite[Theorem~5.9]{ElKGom:inf}, \cite[Theorem~4.15]{GomVer:fir} (cf.\ \cite[Theorem~1.3]{CaeDeG:col}).

  \subsection{Comparison with the results of \cite{NusWam:non}}
 To a given Hopf algebra $H$ and a right $H$-comodule algebra $A$, one associates an $A$-coring $\cC : = A\ot_k H$ with the obvious left $A$-action, diagonal right $A$-action and the coproduct and counit $A\ot_k\Delta_H$ and $A\ot_k\eps_H$ (cf.\ \cite[33.2]{BrzWis:cor}). The category of $(A,H)$-Hopf modules is  isomorphic to the category of right $\cC$-comodules (cf.\ \cite[Proposition~2.2]{Brz:str}). Since $A$ is a right $H$-comodule algebra, it is an $(A,H)$-Hopf module, hence a $\cC$-comodule. The endomorphism ring $\Rend \cC A$ can be identified with the subalgebra of coinvariants $B:= \{b\in A\; |\; \roA(b) = b\ot 1\}$. $A$ is a Galois $\cC$-comodule if and only if $B\subseteq A$ is a Hopf-Galois $H$-extension. With these identifications in mind, the definition of the descent cohomology in \cite{NusWam:non} is a special case of Definition~\ref{def.desc}, \cite[Lemma~2.2]{NusWam:non} can be derived from Lemma~\ref{lemma.action}, \cite[Theorem~2.6]{NusWam:non} follows by Theorem~\ref{thm.desc.twist}, while \cite[Proposition~2.8]{NusWam:non} is a special case of Proposition~\ref{prop.tor}. 

  \subsection{Galois cohomology for partial group actions} As the theory of corings covers all known examples of Hopf-type modules, such as Yetter-Drinfeld modules, Doi-Koppinen Hopf modules, entwined and weak entwined modules, the results of the present note can be easily applied to all these special cases (cf.\ \cite{BrzWis:cor} for more details). One of the special cases of the descent cohomology for corings
   is that of partial Galois actions for non-commutative rings studied in \cite{CaeDeG:par} (cf.\  \cite{DokFer:par}).
    In this section we propose to use the descent cohomology to define the Galois cohomology for partial group actions, thus extending the classical Galois cohomology \cite[Section~I.5]{Ser:Gal}.
   
 Take a finite group $G$ and an algebra $A$. To any element $\sigma\in G$ associate a central idempotent $e_\sigma\in A$ and an isomorphism of ideals $\alpha_\sigma : Ae_{\sigma^{-1}}\to Ae_\sigma$. Following \cite{CaeDeG:par} we say that the collection $(e_\sigma, \alpha_\sigma)_{\sigma\in G}$ is an {\em idempotent partial action}  of $G$ on $A$   if $Ae_1 =A$, $\alpha_1 = A$ and, for all $a\in A$, $\sigma,\tau\in G$,
 $$
 \alpha_\sigma(\alpha_\tau(ae_{\tau^{-1}})e_{\sigma^{-1}}) = \alpha_{\sigma\tau}(ae_{(\sigma\tau)^{-1}})e_\sigma, 
 $$
 (see \cite[Definition~1.1]{DokExe:ass} for the most general definition of a partial group action). Given an idempotent partial action $(e_\sigma, \alpha_\sigma)_{\sigma\in G}$ of $G$ on  $A$, define the invariant subalgebra of $A$,
 $$
 A^G := \{ a\in A\; |\; \forall \sigma\in G, \; \alpha_\sigma(ae_{\sigma^{-1}}) = ae_\sigma\}.
 $$
 The extension $A^G\subseteq A$ is said to be {\em $G$-Galois} if and only if the map
 $$
 A\ot_{A^G}A\to \bigoplus_{\sigma\in G} Ae_\sigma, \qquad a\otimes a'\mapsto \sum_{\sigma\in G} a\alpha_\sigma(a'e_{\sigma^{-1}})v_\sigma
 $$
 is bijective. Here $v_\sigma$ denotes the element of $\bigoplus_{\sigma\in G} Ae_\sigma$ which is equal to $e_\sigma$ at position $\sigma$ and to zero elsewhere. For example, if there exists a convolution invertible right colinear map 
 $k(G)\to A$, where $k(G)$ is the Hopf algebra of functions on $G$, then 
 $A^G\subseteq A$ is a $G$-Galois extension known as a {\em cleft extension} 
 \cite[Section~6.5]{BohVer:Mor}.
 
 By  \cite[Proposition~2.2]{CaeDeG:par} (or as a matter of definition), $(e_\sigma, \alpha_\sigma)_{\sigma\in G}$ is an idempotent partial action of $G$ on $A$ if and only if $\cC := \bigoplus_{\sigma\in G} Ae_\sigma$ is an $A$-coring with the following $A$-actions, coproduct and counit:
 $$
 a(a'v_\sigma)a'' = aa'\alpha_\sigma(a''e_{\sigma^{-1}})v_\sigma, \quad \DC (av_\sigma) = \sum_{\tau\in G} av_\tau\ot_Av _{\tau^{-1}\sigma}, \quad \eC(av_\sigma) = a\delta_{\sigma ,1},
 $$
 for all $a,a',a''\in A$ and $\sigma\in G$. Furthermore, the extension $A^G\subseteq A$ is   $G$-Galois if and only if $\bigoplus_{\sigma\in G} Ae_\sigma$ is a Galois coring (with respect to the grouplike element $\sum_{\sigma\in G} v_\sigma$). 
Following \cite[Definition~2.4]{CaeDeG:par}, a right $A$-module $M$ together with right $A$-module maps $(\varrho_\sigma: M\to M_\sigma)_{\sigma\in G}$ is called a {\em partial Galois descent datum}, provided $\varrho_1 = M$ and each of the $\varrho_\sigma$ restricted to $Me_{\sigma^{-1}}$ is an isomorphism. In view of \cite[Proposition~2.5]{CaeDeG:par}, partial Galois descent data on a right $A$-module $M$
 are in bijective correspondence with descent cycles $Z^1(\bigoplus_{\sigma\in G} Ae_\sigma,M)$. In particular, for any right $A^G$-module $N$, there is a partial Galois descent datum 
 $$
 M=N\ot_{A^G}A, \quad  \varrho_\sigma: N\ot_{A^G}A\to N\ot_{A^G}Ae_\sigma, \quad
 n\ot a\mapsto n\ot \alpha_\sigma(ae_{\sigma^{-1}})e_\sigma.$$
 This is simply the right $\cC$-coaction induced from the coaction on $A$ given by the grouplike element $\sum_{\sigma\in G} v_\sigma$.
  A {\em partial Galois cohomology} of the $G$-Galois extension $A^G\subseteq A$ with values in the automorphism group $\raut AM$ of the right $A$-module $M=N\ot_{A^G}A$ is defined as 
 $$
 H^i (G, \raut A M) := D^i (\bigoplus_{\sigma\in G} Ae_\sigma, M), \qquad i=0,1.
 $$
By Theorem~\ref{thm.desc.twist}, if $A$ is a faithfully flat $A^G$-module, then $H^i (G, \raut A M)$ describes equivalence classes of $A$-twisted forms of $N$.

 As an example  take $N = A^G$. In this case $M=A$ and the automorphism group $\raut AA$ can be identified with the group of units $\cG(A)$ and $\Raut \cC A$ can be identified with $\cG(A^G)$. Furthermore, there is a bijective correspondence between right coactions of $\bigoplus_{\sigma\in G} Ae_\sigma$ on $A$ and the set of grouplike elements in $\bigoplus_{\sigma\in G} Ae_\sigma$ (cf.\ \cite[Lemma~5.1]{Brz:str}). In this way we obtain
$$
H^0(G, \cG(A)) = \cG(A^G), \qquad H^1(G, \cG(A)) = \cG\left(\bigoplus_{\sigma\in G} Ae_\sigma\right)/\cG(A),
$$
where $\cG(A)$ acts on $A$ from the right by conjugation as in Section~\ref{sec.not}.
\section*{Acknowledgement} 
I would like to thank Gabi B\"ohm for discussions and for invitation and warm hospitality in Budapest, where the first version of this note was completed.


\begin{thebibliography}{Bibliography}{}
 \bibitem{BohVer:Mor} G.\ B\"ohm and J.\ Vercruysse, {\em Morita theory for coring extensions and cleft bicomodules},  Adv.\ Math.\ (online) doi:10.1016/j.aim.2006.05.010  (2006).
\bibitem{Brz:str} T.\ Brzezi\'nski,
{\it The structure of corings. Induction functors,
Maschke-type theorem, and Frobenius and Galois-type properties,}
Algebr.\ Represent.\ Theory, 5 (2002), 389--410.
\bibitem{BrzWis:cor} T.\ Brzezi\'nski and R.\ Wisbauer, {\em Corings and
   Comodules}. Cambridge University Press, Cambridge (2003).
\bibitem{Cae:Gal} S.\ Caenepeel,  {\em Galois corings from the descent theory point of view.}  In: {\em Galois theory, Hopf algebras, and semiabelian categories},  Fields Inst.\ Commun., 43, Amer.\ Math.\ Soc., Providence, RI, 2004, pp.\ 163--186.
\bibitem{CaeDeG:par} S.\ Caenepeel and E.\ De Groot, {\em Galois corings applied to partial Galois theory}, in: Proceedings of the International Conference on Mathematics and Applications, ICMA 2004, S.L. Kalla and M.M. Chawla (eds.), Kuwait University, Kuwait 2005. 
    \bibitem{CaeDeG:com} S.\ Caenepeel,
    E.\ De Groot and J.\ Vercruysse, {\em Galois theory for comatrix corings: Descent theory, Morita
    theory, Frobenius and separability properties},  Trans.\ Amer.\ Math.\ Soc., 359 (2007), 185-226.
       \bibitem{CaeDeG:col} S.\ Caenepeel,
    E.\ De Groot and J.\ Vercruysse, {\em Constructing infinite comatrix corings from colimits},  Preprint
    arXiv:math.RA/0511609 (2005), to appear in Appl.\ Categorical Structures.
\bibitem{DokExe:ass} M.\ Dokuchaev and R.\ Exel, {\em Associativity of crossed products by partial actions, enveloping actions and partial representations,} Trans.\ Amer.\ Math.\ Soc.\ 357 (2005), 1931--1952.
\bibitem{DokFer:par} M.\ Dokuchaev, M.\ Ferrero and A.\ Pacques, {\em Partial actions and Galois theory},  J.\ Pure Appl.\ Alg.\  (online) doi:10.1016/j.jpaa.2005.11.009 (2005).
  \bibitem{ElKGom:com} L.\ El Kaoutit and   J.\ G\'omez-Torrecillas, {\em
 Comatrix corings: Galois corings, descent theory, and a structure
theorem for cosemisimple corings.}  Math.\ Z.\ 244 (2003), 887--906.
\bibitem{ElKGom:inf} L.\ El Kaoutit and  J.\ G\'omez-Torrecillas, {\em Infinite comatrix corings.}  Int.\ Math.\ Res.\ Not.\  2004:39 (2004), 2017--2037.
\bibitem{GomVer:fir} J.\ G\'omez-Torrecillas and J.\ Vercruysse,
     {\em Comatrix corings and Galois comodules over firm rings},
     Preprint arXiv:math.RA/0509106 (2005)
\bibitem{NusWam:non}P.\ Nuss and M.\ Wambst, {\em Non-Abelian Hopf cohomology}, Preprint arXiv: math.KT/0511712, (2005)
\bibitem{Ser:Gal} J.-P.\ Serre, {\em Galois Cohomology,} Springer, Berlin (1997)
       \bibitem{Wis:gal}
  R.\ Wisbauer, {\em On Galois comodules,}  Comm.\ Algebra, 34 (2006), 2683--2711.
\end{thebibliography}
\end{document}